\documentclass[12pt]{amsart}
 \usepackage{graphicx}
\usepackage{amssymb}
\usepackage{amsmath}
\usepackage{amscd}
\usepackage{citesort}

\newtheorem{theorem}{Theorem}[section]
\newtheorem{lemma}[theorem]{Lemma}
\newtheorem{e-proposition}[theorem]{Proposition}
\newtheorem{corollary}[theorem]{Corollary}
\newtheorem{e-definition}[theorem]{Definition\rm}
\newtheorem{remark}{\it Remark\/}

\theoremstyle{remark}
\newtheorem{example}[theorem]{{Example}}

\newcommand{\cal}[1]{\mathcal{#1}}
\newcommand{\bez}{\setminus}
\newcommand{\podz}{\subseteq}
\newcommand{\fii}{\varphi}
\newcommand{\sig}{\sigma}
\newcommand{\eps}{\varepsilon}

\newcommand{\fal}[1]{\widetilde{#1}}

\newcommand{\kre}[1]{\overline{#1}}
\newcommand{\gen}[1]{\langle #1 \rangle}
\newcommand{\map}[3]{#1\colon #2\to #3}
\newcommand{\field}[1]{\mathbb{#1}}
\newcommand{\zz}{\field{Z}}
\newcommand{\rr}{\field{R}}

\newcommand{\st}{\;|\;}
\newcommand{\M}[1]{{\cal{M}}(#1)}
\newcommand{\iM}[1]{i_*({\cal{M}}(#1))}

\newcommand{\lst}[2]{{#1}_1,\dotsc,{#1}_{#2}}
\newcommand{\Mob}{M\"{o}bius strip}
\newcommand{\coset}{/\!/}
\newenvironment{tw}{\begin{theorem}}{\end{theorem}}
\newenvironment{lem}{\begin{lemma}}{\end{lemma}}
\newenvironment{prop}{\begin{e-proposition}}{\end{e-proposition}}
\newenvironment{df}{\begin{e-definition}}{\end{e-definition}}
\newenvironment{wn}{\begin{corollary}}{\end{corollary}}
\newenvironment{uw}{\begin{remark}}{\end{remark}}
\newenvironment{ex}{\begin{example}}{\end{example}}
\DeclareMathOperator{\supp}{supp}%

\begin{document}
\title[Commensurability of geometric\ldots]{Commensurability of geometric
subgroups of mapping class groups}
\author{Micha\l\ Stukow}
 \email{trojkat@math.univ.gda.pl}
 \thanks{Supported by the Foundation of Polish Science (FNP)}

\address{Institute of Mathematics, University of Gda\'nsk, Wita Stwosza 57,
80-952 Gda\'nsk, Poland }

\begin{abstract}
Let $M$ be a surface (possibly nonorientable) with punctures and/or boundary components.
The paper is a study of ``geometric subgroups'' of the mapping class group of $M$, that is subgroups corresponding to inclusions of subsurfaces (possibly disconnected). We characterise the subsurfaces which lead to virtually abelian geometric subgroups. We provide algebraic and geometric conditions under which two geometric subgroups are commensurable. We also describe the commensurator of a geometric subgroup in terms of the stabiliser of the underlying subsurface. Finally, we show some applications
of our analysis to the theory of irreducible unitary representations of mapping class groups.
\end{abstract}

\maketitle

\section{Introduction}
Let $M_{g,r}^s$ be a smooth, compact, connected surface of genus $g$ with $s$ punctures and $r$ boundary components (we will call them \emph{holes}). If $r$ and/or $s$ is zero then we omit it from notation and if we do not want to emphasise the numbers $g,r,s$, we simply write $M$ for a surface $M_{g,r}^s$. If $r=0$, we call $M$ a \emph{closed surface}. For the sake of notational convenience we will use the convention that nonorientable surfaces have negative genus, hence $M_{-g}$ is a connected sum of $g$ projective planes, for $g\geq 1$.

If $M$ is a nonorientable surface, define the \emph{mapping class group} $\M{M}$ of $M$ to be the group of isotopy classes of diffeomorphisms of $M$, where we assume that both diffeomorphisms and their isotopies fix the set of punctures and are the identity on the boundary of $M$. The mapping class group of an orientable surface is defined analogously, but we consider only orientation preserving maps. In order to simplify some statements, we define $\M{\emptyset}$ to be the trivial group.

Recall that two subgroups $H_1$ and $H_2$ of a group $G$ are \emph{commensurable} if $H_1\cap H_2$ is of finite index in both $H_1$ and $H_2$. The \emph{commensurator} of $H\leq G$ is defined to be
\[{\rm Comm}(H)=\{g\in G\st \text{$H$ and $gHg^{-1}$ are commensurable}\}. \]

It is not hard to check that if $H_1$ and $H_2$ are commensurable subgroups of $G$, then
\[{\rm Comm}(H_1)={\rm Comm}(H_2).\]
In particular, commensurator is invariant under passing to a finite index subgroup.

The main goal of this paper is to study a family of \emph{geometric subgroups} of $\M{M}$, that is the subgroups of the form $i_*(\M{N})$, where $i_*$ is a homomorphism induced by the inclusion $\map{i}{N}{M}$. To be more precise, $i_*$ is a map that extends a diffeomorphism of a subsurface $N$ of $M$ to a diffeomorphism of $M$. In particular, we
\begin{itemize}
 \item describe the kernel of $i_*$ -- Theorem \ref{tw:injectivity};
 \item describe subsurfaces which lead to virtually abelian geometric subgroups -- Theorem \ref{tw:VA};
 \item give an algebraic and geometric characterisation of geometric subgroups that are commensurable -- Theorems \ref{tw:comm:va}, \ref{tw:GeomtrChar:pre} and \ref{tw:GeomtrChar};
 \item relate the commensurator of a geometric subgroup with the stabiliser of the corresponding subsurface -- Theorems \ref{tw:decr:comm1} and \ref{tw:decr:comm2}.
\end{itemize}
Finally, in section \ref{sec:rep} we provide some straightforward applications of the above results to the theory of unitary representations of mapping class groups -- cf Corollary \ref{cor:induc:repr}.

Some of our results were previously known in the case of a connected subsurface of an orientable surface \cite{RolPar,ParRep}. The novelty of our work is that we allow the subsurfaces to be disconnected and not necessarily injective (i.e. $\map{i_*}{\M{N}}{\M{M}}$ does not need to be injective). The main motivation for this general notion of a geometric subgroup was to include the very important family of subgroups of $\M{M}$, namely the stabilisers of simplexes in the complex of curves on $M$ -- cf Example \ref{ex:complex1}. Moreover, we do not require $M$ to be orientable. The extension to the nonorientable case is possible by the recent results obtained in~\cite{Stukow_twist}.

\section{Preliminaries}
\subsection{Definitions}
By a \emph{circle} in $M$ we mean an unoriented simple closed curve in the interior of $M$ which is disjoint from the set of punctures. Usually we identify a circle with its image. If $a_1$ and $a_2$ are isotopic, we write $a_1\simeq a_2$.
Moreover, as in the case of diffeomorphisms, we will use the same letter for a circle and its isotopy class. By a \emph{boundary circle} we mean a circle parallel to a boundary component of $M$.

According to wether a regular neighbourhood of a circle is an annulus or a \Mob, we call the circle \emph{two-sided} or \emph{one-sided} respectively. We say that a circle is \emph{essential} if it does not bound a disk disjoint from the set of punctures, and \emph{generic} if it bounds neither a disk with fewer than two punctures nor a \Mob\ disjoint from the set of punctures. Notice that the surface $M_{g,r}^s$ admits a generic two-sided circle if and only if $M$ is not $M_{0,r}^s $ with $2r+s\leq 3$ nor $ M_{-1,r}^s$ with $2r+s\leq 2$.

Let $a$ be a two-sided circle. By definition, a regular neighbourhood $S_a$ of $a$ is an
annulus, so if we fix one of its two possible orientations, we can define the
\emph{right Dehn twist} $t_a$ about $a$ in the usual way. We emphasise that since we are
dealing with nonorientable surfaces, there is no canonical way to choose the orientation
of $S_a$. Therefore by a twist about $a$ we always mean one of the two possible twists about
$a$ (the second one is then its inverse). By a \emph{boundary twist} we mean a twist
about a circle parallel to a boundary component. It is known that if $a$ is not
generic then the Dehn twist $t_a$ is trivial. In particular a Dehn twist about the
boundary of a \Mob\ is trivial -- see Theorem 3.4 of \cite{Epstein}.

If $z_1$ and $z_2$ are two punctures in a surface $M$ then there exists their common neighbourhood which is a disk. Hence we can define an \emph{elementary braid} on $z_1$ and $z_2$. It is known that the mapping class group of an orientable surface is finitely generated by Dehn twists and elementary braids \cite{Lick1,Lick2,ParLab,Bir-Punct}.

Other important examples of diffeomorphisms of a nonorientable surface are the
\emph{crosscap slide} and the \emph{puncture slide}. They are defined as a slide of a
crosscap and of a puncture respectively, along a loop. It is known that the mapping class group of a nonorientable surface is finitely generated by Dehn twists, elementary braids, puncture slides and crosscap slides \cite{Lick3,Chil,Kork-non,Stukow_SurBg}.
\subsection{Examples}
It is well known that the mapping class group of $M$ is trivial if and only if $M=M_{0,r}^s$ with $r,s\in \{0,1\}$ or $M=M_{-1,r}^s$ with $s=0$ and $r\leq 1$.

The mapping class group of an annulus or an annulus with a puncture is generated by boundary twists and is isomorphic to $\zz$ or $\zz\times\zz$ respectively. As for less trivial examples,  the mapping class group of a torus or torus with one puncture is generated by twists about meridian and longitude and isomorphic to ${\rm SL}(2,\zz)$. Another nontrivial example is the mapping class group of a disk with $n$ punctures which is isomorphic to the braid group on $n$ strings.

As for nonorientable surfaces, the mapping class group of a projective plane with one or two punctures is generated by puncture slides and is isomorphic to $\zz_2$ or the dihedral group $D_4$ (of order $8$) respectively -- see Corollary 4.6 of \cite{Kork-non}. The mapping class group of a Klein bottle is generated by a twist and a crosscap slide \cite{Lick3}, and is isomorphic to $\zz_2\times\zz_2$. The description of mapping class groups of a Klein bottle with one puncture and a Klein bottle with one hole can be found in the appendix to \cite{Stukow_twist}. In particular, we will use the following proposition.
\begin{prop}\label{prop:but:kl}
Let $M$ be a Klein bottle with one hole. Then
\begin{enumerate}
\item there are exactly two isotopy
classes of generic two-sided circles in $M$, namely the isotopy classes of a boundary circle $b$ and of a nonseparating two-sided circle $a$ in $M$ (see Figure~\ref{RN02}); 
\begin{figure}[h]
\includegraphics{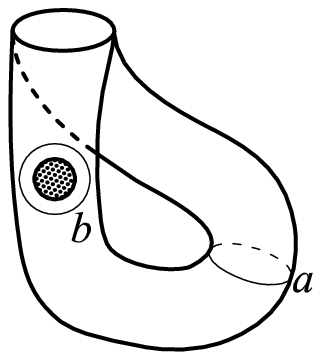}
\caption{Generic two-sided circles on a surface $M_{-2,1}$.} \label{RN02}
\end{figure}
\item the group generated by
Dehn twists $t_a$ and $t_b$ is an index two subgroup of $\M{M}$ and is isomorphic to $\zz\times\zz$.
\end{enumerate}
\end{prop}
\begin{proof}
 The idea of the proof of assertion (1) is very simple, first one observes that if $c$ and $d$ are generic two-sided circles of
 the same separability (i.e. both are separating or nonseparating) then there exists a
 diffeomorphism $\map{h}{M}{M}$ such that $h(c)\simeq d$. Then from the structure of the mapping
 class group of $M$ (cf Theorem A.7 of \cite{Stukow_twist}) one concludes that for every
 $h\in{\cal{M}}(M)$, $h(a)\simeq a$ and $h(b)\simeq b$. We omit the details, refereing the
 reader to the fully analogous proof of Proposition A.3 in \cite{Stukow_twist}.

 Assertion (2) is a consequence of Theorem A.7 of \cite{Stukow_twist}.
\end{proof}
\subsection{Pantalon \& skirt decompositions}
\begin{figure}[h]
\includegraphics{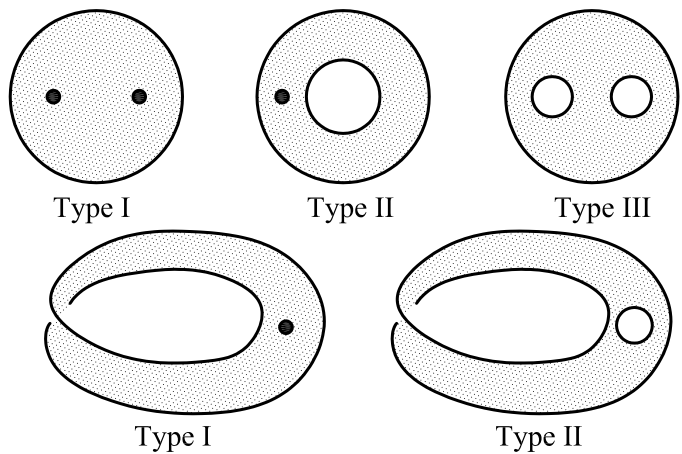}
\caption{Different types of pantalons and skirts.} \label{fig:pant}
\end{figure}

Following \cite{RolPar}, we call the surfaces $M_{0,1}^2$, $M_{0,2}^1$ and $M_{0,3}$ pantalons of type I, II and III respectively (cf Figure \ref{fig:pant}).
We say that a collection of different two-sided circles $\lst{a}{n}$ on $M$ define a \emph{pantalon decomposition} of $M$ if each connected component of $\kre{M\bez\bigcup_{i=1}^na_i}$ is a pantalon. It is known \cite{RolPar} that an orientable surface has a pantalon decomposition if and only if $2g+r+s> 2$ and $M\neq M_0^3$.

As observed in Section 5 of \cite{Stukow_twist}, one needs to add two more ``pieces'' in order to decompose nonorientable surfaces, namely a \Mob\ with one puncture and a \Mob\ with an open disk removed. We call these surfaces a \emph{skirt of type I and II}, respectively (cf Figure \ref{fig:pant}). A decomposition of a surface into pantalons and skirts is called a \emph{P-S decomposition}. The nonorientable surface $M_{g,r}^s$ has a P-S decomposition if and only if $r+s-g>2$ and $M\neq M_{-1}^2$.

The mapping class groups of pantalons of type II, III and of a skirt of type II are generated by boundary twists and isomorphic to $\zz^2$, $\zz^3$ and $\zz^2$ respectively. The mapping class groups of a pantalon of type I and of a skirt of type I are isomorphic to $\zz$ and generated by an elementary braid and a puncture slide along a one-sided loop respectively.
\subsection{Properties of circles}
For any two circles $a$ and $b$, we define their \emph{geometric intersection number}:
\[I(a,b)={\rm inf}\{|a'\cap b|\,: a'\simeq a\}. \]
In particular, if $a$ is a two-sided circle and $a\simeq b$ then $I(a,b)=0$.
\begin{prop}[Propositions 3.8 of \cite{RolPar} and 4.4 of \cite{Stukow_twist}]
\label{prop:4.4} Suppose $M\neq M_{-2}$ and let $\lst{a}{n}$ be generic, pairwise nonisotopic, pairwise disjoint and two-sided
circles on $M$. Then the homomorphism \[\map{h}{\zz^n}{{\cal{M}}(M)}\]
defined by \[h(\lst{\alpha}{n})=t_{a_1}^{\alpha_1}\cdots t_{a_n}^{\alpha_n}\] is injective. \qed
\end{prop}
\begin{prop}[Propositions 3.7 of \cite{RolPar} and 4.7 of \cite{Stukow_twist}]\label{prop:4.7} Let $a$ and $b$ be generic two-sided circles in $M$. If
$\alpha$ and $\beta$ are nonzero integers such that $t_a^\alpha t_b^\beta=t_b^\beta t_a^\alpha$, then $I(a,b)=0$.
\qed
\end{prop}

Although the following proposition is stated in \cite{RolPar} only
for orientable surfaces, the proof applies to the nonorientable case verbatim.

\begin{prop}[Proposition 6.2 of \cite{RolPar}]\label{prop_6_2}
Suppose $\lst{a}{n}$ are essential two-sided circles which are pairwise disjoint. Let $b$ be an
essential two-sided circle such that $I(a_i,b)=0$ for all $i=1,\ldots,n$. Then there exists a
two-sided circle $c$ isotopic to $b$ such that $a_i\cap c=\emptyset$ for all $i=1,\ldots,n$.  \qed
\end{prop}

%
%

\begin{prop}\label{prop_6_10}
Let ${\cal{A}}=\{\lst{a}{n}\}$ and ${\cal{B}}=\{\lst{b}{n}\}$ be two collections of pairwise disjoint,
essential and two-sided circles on $M$ such that $a_i$ is isotopic to $b_i$, for each $i=1,\ldots,
n$. Then there exists an isotopy $\map{h_t}{M}{M}$ such that $h_0$ is an identity and
$h_1({\cal{A}})={\cal{B}}$.
\end{prop}
\begin{proof}
We will use induction on the number $m$ of isotopy classes of circles that have more than one
 representant in $\cal{A}$. For $m=0$ the statement
follows from Proposition 3.10 of \cite{RolPar} (although the statement in \cite{RolPar} is only
for orientable surfaces, the proof applies to the nonorientable case verbatim).

Now assume that $m\geq1$ and up to the permutation of elements of ${\cal{A}}$ and ${\cal{B}}$ we can
assume that $a_1\simeq a_2\simeq\ldots\simeq a_k$ for $2\leq k\leq n$ and
$a_1\not\simeq a_j$ for $j> k$. Applying the inductive hypothesis to the collections ${\cal
A}'=\{a_k,a_{k+1},\ldots,a_n\}$ and ${\cal B}'=\{b_k,b_{k+1},\ldots,b_n\}$ we obtain an isotopy
$\map{h_t'}{M}{M}$ such that $h'_0$ is the identity, $h'_1({\cal{A}}')={\cal{B}}'$ and
$h'_1(a_k)=b_k$. Clearly the isotopies $a_1\simeq a_2\simeq\ldots\simeq a_k$
provide a family of annuli between the circles $\{a_1,\ldots,a_k\}$. Let $A$ be a maximal one
with respect to inclusion. 
Hence $\{a_1,\ldots,a_k\}\podz A$ and $\partial A \podz \{a_1,\ldots,a_k\}$. If we define $B$ in a
similar manner but with respect to the circles $\{b_1 \ldots,b_k \}$, then $h'_1(A)$ and $B$, as regular neighbourhoods of $h_1'(a_k)=b_k$,  are
isotopic by an isotopy $\map{h_t''}{M}{M}$. Moreover, since both $h'_1(A)$ and $B$ are
disjoint from $b_{k+1},\ldots,b_n$, we can assume that $h''(b_i)=b_i$ for $i=k+1,\ldots,n$. Finally,
let $\map{h_t'''}{M}{M}$ be an isotopy of $B=h_1''h_1'(A)$ which transforms
$h''_1(h'_1(\{\lst{a}{k}\}))$ onto $\{\lst{b}{k}\}$. Then the composition
$h_t(x)=h'''_t(h''_t(h'_t(x)))$ is a required isotopy between $\cal{A}$ and $\cal{B}$.
\end{proof}
\begin{lem}\label{lem:zaw:but:kle}
Let $a_0$ and $a_1$ be two disjoint, generic, nonisotopic and separating circles on a surface $M$ such that one of the connected components of $\kre{M\bez a_i}$ is a Klein bottle $K_i$ with one hole for $i=0,1$. Then $K_0\cap K_1=\emptyset$.
\end{lem}
\begin{proof}
Suppose that $K_0\cap K_1\neq \emptyset$. If $K_0\podz K_1$ or $K_1\podz K_0$ then one easily concludes that $a_0\simeq a_1$, hence we can assume that this is not the case. Therefore $a_0$ intersects $K_1$, hence it is contained in the interior of $K_1$. By part (1) of Proposition \ref{prop:but:kl}, this implies that $a_0\simeq a_1$, a contradiction.
\end{proof}



\subsection{Some basic properties of a subgroup index}

For the sake of completeness we review below some basic properties of a subgroup index. Throughout this section we will use the notation $G\coset H$ for the set of left cosets of $H$ in $G$.
\begin{prop}\label{prop:ex:virt}
Let $\map{\fii}{G}{\fii(G)}$ be a group homomorphism and assume that $H$ is a finite index subgroup of $G$. Then $\fii(H)$ is a finite index subgroup of $\fii(G)$. In particular, a homomorphic image of a virtually abelian group is virtually abelian.
\end{prop}
\begin{proof}
Define a map
 $\map{\Phi}{\fii(G)\coset\fii(H)}{G\coset H}$ as follows:
 \[\Phi(\fii(g)\fii(H))=gH, \]
 where we choose one $g\in G$ for each coset in $\fii(G)\coset\fii(H)$. It is straightforward to check that $\Phi$ is ``1-1''.
\end{proof}
\begin{prop}\label{prop:ex:cap}
Let $H$ be a finite index subgroup of a group $G$, and let $K\leq G$ be any subgroup. Then $H\cap K$ has finite index in $K$.
\end{prop}
\begin{proof}
Define a map $\map{\Phi}{K\coset(H\cap K)}{G\coset H}$ as follows:
\[\Phi(g(H\cap K))=gH, \]
where we choose one $g\in G$ for each coset in $K\coset (H\cap K)$. It is straightforward to check that $\Phi$ is ``1-1''.
\end{proof}
\begin{prop}\label{prop:ex:ext:prod}
Let $H_i\leq G_i$  be a finite index subgroup for $i=0,1$. Then $H_0\times H_1$ has finite index in $G_0\times G_1$.
\end{prop}
\begin{proof}
The assertion follows from the well know (and easy to check) formula:
\[[G_0\times G_1:H_0\times H_1]=[G_0:H_0][G_1:H_1] \]
\end{proof}
\begin{prop}\label{prop:ex:semprod}
Let $H_0,H_1,K$ be subgroups of a group $G$ such that $K$ centralises both $H_0$ and $H_1$. If $H_0$ and $H_1$ are commensurable then $H_0K$ and $H_1K$ are also commensurable.
\end{prop}
\begin{proof}
Since $H_0$ and $H_1$ are commensurable, by Proposition \ref{prop:ex:ext:prod}, \[(H_0\cap H_1)\times K\] is a finite index subgroup of both $H_0\times K$ and $H_1\times K$. Moreover, since $K$ centralises $H_0$ and $H_1$, we have homomorphisms $\map{\fii_i}{H_i\times K}{G}$ defined by $\fii_i(h,k)=hk$, for $i=0,1$. Therefore by Proposition \ref{prop:ex:virt}, $(H_0\cap H_1)K$ is a finite index subgroup of both $H_0K$ and $H_1K$. This finishes the proof, since it is straightforward to check that
\[(H_0\cap H_1)K\leq H_0K\cap H_1K \]
\end{proof}
The following example shows that the implication in the statement of Proposition \ref{prop:ex:semprod} can not be in general replaced by an equivalence, even if we assume that $H_iK=H_i\oplus K$ for $i=0,1$.
\begin{ex}\label{ex:semprod:opp}
Let $G=\gen{a,b\st[a,b]}$ be a free abelian group of rank $2$, and let $H_0=\gen{a}$, $H_1=\gen{ab}$, $K=\gen{b}$. Then
\[G=H_0\oplus K=H_1\oplus K\]
but $H_0\cap H_1=1$.
\end{ex}

\section{Subsurfaces and injectivity}
Following \cite{RolPar}, define an \emph{exterior cylinder} $E$ of a subsurface
${N\subset M}$ to be any component of $\kre{M\bez N}$ that is an annulus with both
boundary circles $a$ and $b$ in $N$. In what follows, we will always assume that
the orientations of regular neighbourhoods of $a$ and $b$ agree with the chosen
orientation of $E$. In other words, twists $t_{a}$ and $t_{b}$ are equal in
${\cal{M}}(E)$.
\begin{df}
Let $N\neq M$ be a closed subsurface of $M$ (not necessarily connected). We call $N$ an \emph{essential
subsurface} if the following conditions are satisfied
\begin{enumerate}
 \item every boundary component of $N$ is generic in $N$ and $N$ is disjoint from $\partial M$;
 \item for every connected component $C$ of $N$ which is an annulus, the meridian of $C$ is not
 isotopic in $M$ to a boundary component of $N\bez C$;
 \item no component of $\kre{M\bez N}$ is a disk.
\end{enumerate}
\end{df}
The first of the above conditions means that $N$ has no components with trivial mapping class
group, that is components homeomorphic to a disk with less than two punctures or a \Mob. The
second one implies that the generator of the mapping class group of $C$ does not belong to
${\cal M}(N\bez C)$. Therefore, from a mapping class group point of view, these two assumptions
are quite natural and it turns out that they greatly simplify some arguments. The third
condition is technical and it implies the following proposition.
\begin{prop}[Proposition 3.5 of \cite{RolPar}]\label{prop:3_5}
Let $N$ be an essential subsurface of $M$ and let $a,b$ be essential two-sided circles in $N$ such that $a$ in not a boundary circle of an exterior cylinder to $N$. Then $a$ and $b$ are isotopic in $M$ if and only if they are isotopic in $N$.\qed
\end{prop}
As an immediate corollary we obtain:
\begin{prop}
Let $a,b$ be essential circles in an essential subsurface $N$ of $M$, and let $I_N(a,b)$ denote the geometric intersection number of $a$ and $b$ treated as circles in $N$. Then $I_N(a,b)=I(a,b)$.\qed
\end{prop}
Keeping in mind the above propositions we will often abuse notation by identifying isotopy in $N$ with isotopy in $M$ and $I_N(a,b)$ with $I(a,b)$.


In some of our applications we will need to impose some further conditions on subsurfaces.
\begin{df}
We call an essential subsurface $N\subset M$ \emph{generic} if every boundary component of $N$ is generic in $M$ (hence does not bound a disk with less than 2 punctures nor a \Mob).
\end{df}
\begin{ex}\label{ex:complex1}
Let $X_n(M)$ be the set of $n+1$-tuples $(a_0,\ldots,a_{n})$ of one-sided and generic two-sided disjoint circles in $M$, such that $a_i$ is not a boundary circle of $M$ and $a_i$ is not isotopic to $a_j$ for $i\neq j$. We say that two elements $(a_0,\ldots,a_n)$ and $(b_0,\ldots,b_n)$ of $X_n(M)$ are equivalent if there exists a permutation $\map{s}{\{0,\ldots,n\}}{\{0,\ldots,n\}}$ such that $a_i$ is isotopic to $b_{s(i)}$ for $i=0,\ldots, n$.

Define the \emph{complex of curves} $C(M)$ to be an abstract simplicial complex such that the $n$-simplexes in $C(M)$ are the equivalence classes of elements of $X_n(M)$ with respect to the equivalence relation defined above. Clearly we can think of elements of $C(M)$ as sets of isotopy classes of pairwise disjoint circles.

Every element $\sig\in C(M)$ provides a natural example of an essential subsurface $M_\sig$ of $M$, namely $M_\sig$ is the complement of a regular neighbourhood of $\sig\cup\partial M$. Observe that $M_\sig$ is a generic subsurface if and only if $\sig$ does not contain one-sided circles.
\end{ex}


Consider an essential subsurface $N\subset M$ and the homomorphism
$\map{i_*^N}{{\cal{M}}(N)}{{\cal{M}}(M)}$ induced by inclusion. If we do not want to emphasise the subsurface $N$, we simply write $i_*$ for $i_*^N$. The following theorem describes the kernel of $i_*$.


 \begin{tw}\label{tw:injectivity}
Let $N$ be an essential subsurface of $M\neq M_{-2}$ and let $\lst{C}{p}$ be components
of $N$ that are annuli. Denote by $\lst{a}{n}$ the boundary components of
$N\bez{\bigcup_{i=1}^{p}C_i}$ that are not generic in $M$, and let $\lst{c}{m}$ be these
meridians of $\lst{C}{p}$ that are not generic in $M$. Denote also by $b_i$, $b_i'$, for
$i=1,\ldots,k$, the pairs of boundary components of $N\bez{\bigcup_{i=1}^{p}C_i}$ that
bound exterior cylinders. Then the kernel of $i_*^N$ is generated by
$\{t_{a_1},\ldots,t_{a_n},t_{c_1},\ldots,t_{c_m},t_{b_1}^{-1}t_{b'_1},\ldots,t_{b_k}^{-1}t_{b'_k}\}$. Moreover, this kernel is isomorphic to $\zz^{n+m+k}$. \qed
\end{tw}
\begin{proof}
The proof
is very similar to the proof of Theorem 4.1 of \cite{RolPar}, so we only give the main idea,
refereing the reader to \cite{RolPar}.

Since $N$ is essential, $N\bez \bigcup_{i=1}^{p} C_i$ has a P-S decomposition. Let $\lst{d}{u}$
be the union of circles defining this decomposition and these meridians of $\lst{C}{p}$ that
are generic in $M$. Let $\lst{e}{w}$ be the boundary components of $N\bez \bigcup_{i=1}^{p} C_i$
different from all of $a_i$ and $b_i,b_i'$ (see Figure \ref{RN01}).
\begin{figure}[h]
\includegraphics{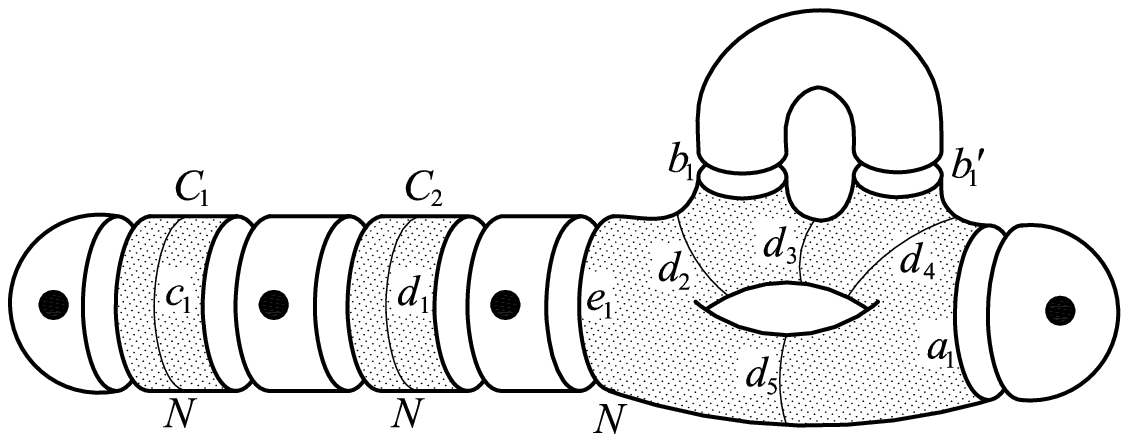}
\caption{Subsurface $N$ -- Theorem \ref{tw:injectivity}.} \label{RN01}
\end{figure}
%

Let $h\in \ker i_*$. Since $h$, as an element of $\M{M}$, is isotopic to the identity in $M$, $h(d_i)\simeq d_i$ in $M$
for $i=1,\ldots,u$. By Proposition \ref{prop:3_5}, $h(d_i)\simeq d_i$ in $N$, hence by
Proposition \ref{prop_6_10}, we can assume that $h$, as an element of $\M{N}$, is the identity on each of $\lst{d}{u}$
and on boundary curves of $N$. Moreover, since isotopies fix punctures, $h$ does not permute them
and does not reverse the local orientations around them. Therefore, by the structure of mapping
class groups of pantalons, skirts and annulus, we conclude that
$$h=t_{a_1}^{\alpha_1}\cdots t_{a_n}^{\alpha_n}t_{b_1}^{\beta_1}t_{b_1'}^{\beta_1'}\cdots
t_{b_k}^{\beta_k}t_{b_k'}^{\beta_k'}t_{c_1}^{\gamma_1}\cdots
t_{c_m}^{\gamma_m}t_{d_1}^{\delta_1}\cdots t_{d_u}^{\delta_u}t_{e_1}^{\eps_1}\cdots
t_{e_w}^{\eps_w},$$ for some integers $\alpha_i,\beta_i,\beta_i',\gamma_i,\delta_i,\eps_i$. Therefore
$$1=i_*(h)=t_{b_1}^{\beta_1+\beta_1'}\cdots
t_{b_k}^{\beta_k+\beta_k'}t_{d_1}^{\delta_1}\cdots t_{d_u}^{\delta_u}t_{e_1}^{\eps_1}\cdots
t_{e_w}^{\eps_w}.$$

By Proposition \ref{prop:4.4},
 $$\beta_1+\beta_1'=\ldots=\beta_k+\beta_k'=\delta_1=\ldots=\delta_u=\eps_1=\ldots=\eps_w=0,$$
 which proves that $\ker
 i_*=\gen{t_{a_1},\ldots,t_{a_n},t_{c_1},\ldots,t_{c_m},t_{b_1}^{-1}t_{b'_1},\ldots,t_{b_k}^{-1}t_{b'_k}}$, and $\ker i_*\cong\zz^{n+m+k}$.
\end{proof}
\begin{df}
If the homomorphism $\map{i_*}{{\cal{M}}(N)}{{\cal{M}}(M)}$ is injective, we call
$N$ an \emph{injective subsurface}. In such a case we usually identify
$i_*({{\cal{M}}(N)})$ with ${{\cal{M}}(N)}$.
\end{df}
An injective subsurface is of course generic, however a generic subsurface can have exterior cylinders.
\begin{wn}
Let $N$ be an essential subsurface of $M$ and $M$ is not a torus. Then $N$ is injective if and only
if no component of $\kre{M\bez N}$ is a disk with less then two
punctures, \Mob\ or an annulus whose boundary components are both boundary components of
$N$. \qed
\end{wn}
Observe that the above corollary can be thought as a generalisation of Proposition
\ref{prop:4.4}. In fact, Proposition \ref{prop:4.4} follows for $N$ being an union of
disjoint regular neighbourhoods of appropriate circles.




\begin{wn}\label{rem:prod}
Let $\lst{U}{p}$ be connected components of an essential subsurface $N$ of $M$. Then the geometric subgroup $i_*({\cal{M}}(N))$
is isomorphic to the quotient of the product $\prod_{i=1}^p
\iM{U_i}$ by the subgroup generated by
$\{t_{b_1}^{-1}t_{b'_1},\ldots,t_{b_k}^{-1}t_{b'_k}\}$,
where $b_1,b_1',\ldots,b_k,b_k'$ are pairs of boundary components of $N$ that
bound exterior cylinders. \qed
\end{wn}

\begin{prop}\label{prop:gencirc}
Let $N$ be an essential subsurface of $M$ and assume that $N$ is not a Klein bottle with one hole.
Let $a$ be a two-sided generic circle in $N$ which is not a boundary circle of
$N$. Then there exists a two-sided circle $b$ in $N$, which is generic in $M$, and  $I(a,b)>0$.
\end{prop}
\begin{proof}
For $M=M_{-2}$ the statement is trivial, hence assume that $M\neq M_{-2}$. By Proposition 3.4 of \cite{RolPar} and Lemma 4.1 of \cite{Stukow_twist} there exists a circle $b$, generic in $N$, such that $I(a,b)>0$. In particular $b$ is not a boundary circle of $N$, hence by Proposition \ref{prop:4.4} and Theorem \ref{tw:injectivity}, the twist $t_b$ is nontrivial in $M$, that is $b$ is generic in $M$.
\end{proof}


\section{Diffeomorphisms of subsurfaces and their isotopies}
The main goal of this section is to prove Proposition \ref{prop:supp:disjoint2} which, roughly speaking, characterises isotopy classes of diffeomorphisms of an essential subsurface $N$ which have representants with support disjoint from $N$. This result will be an essential tool in proving a partial converse to Proposition \ref{prop:ex:semprod} -- see Lemma \ref{lem:comm:red:common}.
\begin{lem}\label{prop:supp:disjoint}
Let $U$ be a sum of some of the connected components of an essential subsurface $N$ of $M$ (see Figure \ref{RN08}).
Suppose that $f\in\M{M}$ is such that $f\in i_*(\M{N})$ and the support of some representant of $f$ is disjoint from $U$. Then 
\[f= gt_{b_1}^{\beta_1}\cdots t_{b_k}^{\beta_k}, \]
where $g\in \iM{N\bez U}$ and each of $\lst{b}{k}$ is a boundary circle of~$U$.
\begin{figure}[h]
\includegraphics{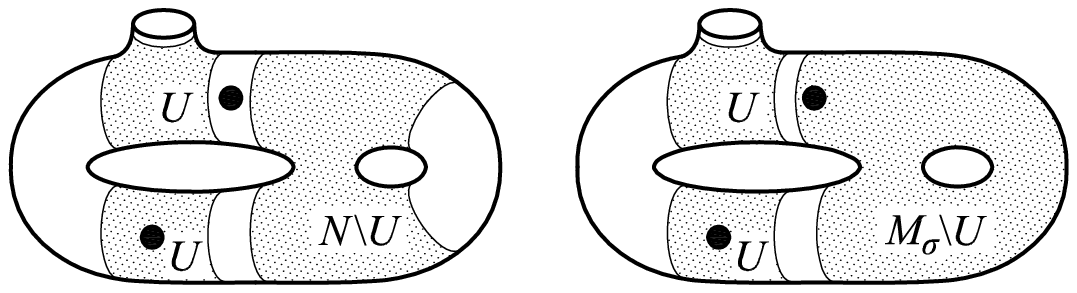}
\caption{Subsurfaces $N$ and $M_\sig$ - Lemma \ref{prop:supp:disjoint}.} \label{RN08}
\end{figure}
\end{lem}
\begin{proof}
Since the proof is quite technical, we will first explain the main idea. We know that the isotopy class of $f$ contains two representants $f_1$ and $f_2$ such that $f_1$ acts trivially outside $N$ and $f_2$ acts trivially on $U$. What we are going to show is that in fact, up to the twists about boundary components of $U$, $f$ can be represented by a diffeomorphism $g$ which shares these two properties, that is the support of $g$ is contained in $N\bez U$. The idea of the construction of $g$ is to glue the action of $f_1$ on $N\bez U$ with the action of $f_2$ on $U$. In order to do this rigorously, we will use Theorem \ref{tw:injectivity}.

Observe that if $U$ contains a component $C$ homeomorphic to an annulus and  $U'=U\bez C$, then under the same assumptions about $f$, it is enough to prove that
\[f\simeq g't_{b_1}^{\beta_1}\cdots t_{b_k}^{\beta_k}, \]
where $g'\in \iM{N\bez U'}$ and $b_1,\ldots,b_k$ are boundary circles of $U'$. In fact, this is an easy consequence of the fact that the mapping class group of $C$ is generated by a boundary twist. Therefore we can assume that $U$ does not contain an annulus as a component. Assume also that $U\neq \emptyset$ and $M\bez U$ is not a collection of disjoint annuli -- in such cases the assertion is obvious.

Let $\sig$ be the set of isotopy classes of the boundary circles of $U$ which are not boundary circles of $M$.
 Define $M_\sig$ to be a complement of $\sig$ in $M$ as in Example \ref{ex:complex1} -- see Figure \ref{RN08}. By an appropriate choice of regular neighbourhoods of elements of $\sig$ we can assume that $U$ is a sum of connected components of $M_\sig$ (we use here the assumption that $U$ does not have annuli as connected components). Moreover, since $M\bez U$ is not a collection of annuli, $M_\sig\neq U$.


Since $f\in \iM{N}$, there exists a diffeomorphism $f_1\in \M{N\cap M_\sig}$ such that $i_*^{M_\sig}(f_1)=f$. On the other hand, $f$ has a representant $f_2$ with support disjoint from $U$. By composing $f$ with some powers of twists about meridians of the components of $M\bez M_\sig$ which are annuli, we can assume that the action of $f_2$ on each such annulus is trivial. Hence we can treat $f_2$ as an element of $\M{M_\sig\bez U}$.
The basic properties of $f_1$ and $f_2$ are as follows:
\begin{enumerate}
\item $i_*^{M_\sigma}(f_1)=i_*^{M_\sigma}(f_2)=f$,
\item the action of $f_2$ on $U$ is trivial,
\item the support of $f_1$ is contained in $N$.
\end{enumerate}
By Theorem \ref{tw:injectivity} and assertion (1) above,
\[f_1\simeq f_2t_{b_1}^{\beta_1}t_{b'_1}^{-\beta_1}\cdots t_{b_k}^{\beta_k}t_{b'_k}^{-\beta_k}u, \]
where $u$ is a product of some powers of boundary twists in $U$ such that $i_*(u)=1$, $\lst{b}{k}$ and $b'_1,\ldots,b'_k$ are boundary circles in $U$ and ${M_\sigma\bez U}$ respectively such that $b_i$ is isotopic to $b'_i$ in $M$, for $i=1,\ldots,k$. By property (2),
we obtain that
\[g\simeq f_1t_{b_1}^{-\beta_1}\cdots t_{b_k}^{-\beta_k}u^{-1}
\simeq f_2t_{b'_1}^{-\beta_1}\cdots t_{b'_k}^{-\beta_k}\] is an element of $\M{M_\sig}$ which acts trivially on $U$.
By Corollary \ref{rem:prod},
\[\M{M_\sig}=\M{M_\sig\bez U}\oplus \M{U},\] hence we can change the action of $g$ on $U$ to the identity without changing its action on $M_\sig\bez U$. In other words, using (3), we can assume that the support of $g$ is contained in $N\bez U$.
Therefore
\[f=i_*(f_1)=i_*(gt_{b_1}^{\beta_1}\cdots t_{b_k}^{\beta_k}u)=gt_{b_1}^{\beta_1}\cdots t_{b_k}^{\beta_k}, \]which completes the proof.
%
%
\end{proof}

\begin{lem} \label{lem:comm:pow:tw}
Let $N$ be an essential subsurface of $M$, and let $\lst{a}{n}$ be pairwise disjoint, generic and two-sided circles on $M$ such that \mbox{$a_i\cap N=\emptyset$} for $i=1,\ldots,n$.
If the product of twists \[t_{a_1}^{\alpha_1} \cdots t_{a_n}^{\alpha_n}\] is an element of $\iM{N}$ for nonzero integers $\alpha_1,\ldots,\alpha_n$, then each of the $\lst{a}{n}$ is isotopic
in $M$ to a boundary component of $N$.
\end{lem}
\begin{proof}
If $M=M_{-2}$ then the assertion follows from the fact that there is only one generic
circle in $M_{-2}$ -- see \cite{Lick3}. Therefore assume that $M\neq M_{-2}$.

If we apply Lemma \ref{prop:supp:disjoint} with $U=N$ and $f=t_{a_1}^{\alpha_1}\cdots t_{a_n}^{\alpha_n}$, we obtain
\[t_{a_1}^{\alpha_1}\cdots t_{a_n}^{\alpha_n}=t_{b_1}^{\beta_1}\cdots
t_{b_k}^{\beta_k}, \]
where each of $b_1,\dots,b_k$ is a boundary circle of $N$ which is generic in $M$.
By Proposition
\ref{prop:4.4}, each $a_i$ is isotopic to some $b_j$.
\end{proof}

\begin{lem}\label{lem:supp:annulus:weg}
Let $C$ be a connected component of an essential subsurface $N$ of $M$ which is an annulus with meridian $c$. Let $K$ be an essential subsurface of $M$ which do not have boundary circles isotopic in $M$ to $c$. Suppose that $f\in\M{M}$ is such that $f\in i_*(\M{N})$ and the support of some representant of $f$ is contained in $K$. Then $f\in\iM{N\bez C}$.
\end{lem}
\begin{proof}
By Lemma \ref{prop:supp:disjoint},
\begin{equation}f\simeq gt_{c}^{\gamma}\label{eq:lemm:supp},\end{equation}
where $g\in \iM{N\bez C}$. 

Since both $f$ and $t_c^\gamma$ act trivially on $\partial K$, we can decompose $g\simeq g_1g_2$, for diffeomorphisms $g_1,g_2\in \iM{N\bez C}$ such that $\supp(g_1)\podz N\bez K$ and $\supp(g_2)\podz K$. Clearly we can assume that $C$ is disjoint from $K$, hence $\fal{K}=K\cup C$ is an essential subsurface of $M$. Since $f,g_2,t_c^\gamma\in \iM{\fal{K}}$, equation \eqref{eq:lemm:supp} implies that $g_1\in\iM{\fal{K}}$. Hence applying Lemma \ref{prop:supp:disjoint} with $f=g_1$, $N=\fal{K}$ and $U={K}$, we obtain
\[g_1\simeq t_{b_1}^{\beta_1}\cdots t_{b_k}^{\beta_k}, \]
where each of $b_1,\ldots,b_k$ is a boundary circle of $\fal{K}$. We claim that none of these circles is isotopic to $c$. In fact, since
\[t_{b_1}^{\beta_1}\cdots t_{b_k}^{\beta_k}\simeq g_1 \in \iM{N\bez C}, \]
 by Lemma \ref{lem:comm:pow:tw}, we obtain that each of $b_1,\ldots, b_k$ is a boundary circle of $N\bez C$.

Therefore
\[t_{b_1}^{\beta_1}\cdots t_{b_k}^{\beta_k}t_c^{\gamma}\simeq g_1t_c^{\gamma}\simeq g_1 g^{-1}f\simeq g_2f\in \iM{K}.\]
Since $b_i\not\simeq c$ for $i=1,\dots,k$, Lemma \ref{lem:comm:pow:tw} implies that $\gamma=0$.
\end{proof}
\begin{prop}\label{prop:supp:disjoint2}
Let $U$ be a sum of some of the connected components of an essential subsurface $N$ of $M$, and let $K$ be an essential subsurface of $M$ such that $K\cap U=\emptyset$ and $K\cup U$ is essential.
Suppose that $f\in\M{M}$ is such that $f\in i_*(\M{N})$ and the support of some representant of $f$ is contained in $K$. Then 
\[f=gt_{b_1}^{\beta_1}\cdots t_{b_k}^{\beta_k}, \]
where $g\in \iM{N\bez U}$ and each of $\lst{b}{k}$ is a boundary circle of both $U$ and $K$.
\end{prop}
\begin{proof}
By Lemma \ref{prop:supp:disjoint},
\[f=gt_{b_1}^{\beta_1}\cdots t_{b_k}^{\beta_k}t_{c_1}^{\gamma_1}\ldots t_{c_m}^{\gamma_m}, \]
where $g\in \iM{N\bez U}$, $b_1,\ldots, b_k$ are the boundary circles of both $U$ and $K$ and $\lst{c}{m}$ are boundary circles of $U$ which are not boundary circles of $K$. Hence it is enough to show that $\gamma_1=\cdots=\gamma_m=0$.

Let $C_i$ be a regular neighbourhood of $c_i$ disjoint from $N\bez U$ for $1 \leq i\leq m$, and let
\[\fal{N}=(N\bez U)\cup\bigcup_{i=1}^m C_j.\] 
For a fixed $1 \leq j\leq m$, applying Lemma \ref{lem:supp:annulus:weg} with $N=\fal{N}$, $C=C_j$, $K=K$ and $f=f$, we obtain that $f\in \iM{\fal{N}\bez C_j}$. Hence, by Corollary \ref{rem:prod}, $\gamma_j=0$.
\end{proof}
The following simple example is an attempt to convince the reader that the above proposition is quite interesting in its own right because it provides some constraints on possible relations in the mapping class group.
\begin{ex}\label{ex:torus_two_holes}
Suppose that $M$ is a torus with two holes and let the circles $a_1,a_2,a_3,b_1,b_2$ be as in Figure \ref{RN07}(i).
\begin{figure}[h]
\includegraphics{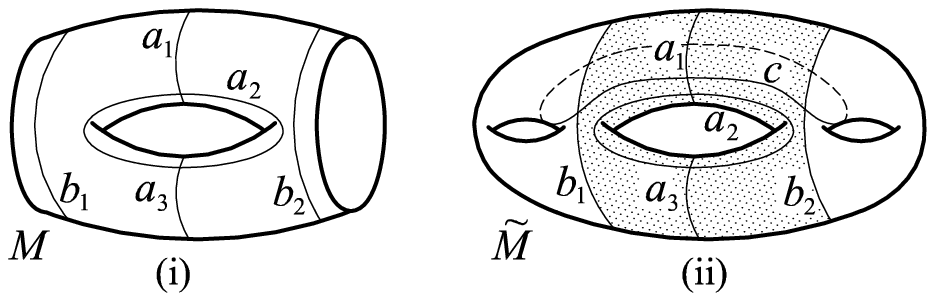}
\caption{Torus with $2$ holes - Example \ref{ex:torus_two_holes}.} \label{RN07}
\end{figure}
Then, with an appropriate choice of orientations of regular neighbourhoods of these circles, there is a ``torus with two holes'' relation in the mapping class group \cite{Kork-non1}, namely:
\[(t_{a_1}t_{a_2}t_{a_3})^4=t_{b_1}t_{b_2}.\]
We claim that Proposition \ref{prop:supp:disjoint2} implies that the above relation can not be ``resolved'' with respect to any of the twists on the left hand side. To be more precise, we concentrate on the twist $t_{a_1}$. We claim that if $\alpha\neq 0$ then $t_{a_1}^\alpha$ is not an element of the group generated by the remaining twists. In fact, otherwise $N$=regular neighbourhood of $a_2\cup a_3\cup b_1\cup b_2$, $K$=regular neigbourhood of $a_1$, $U$=regular neighbourhood of $b_1\cup b_2$ and $f=t_{a_1}^\alpha$ would satisfy the assumptions of Proposition \ref{prop:supp:disjoint2}, hence
\[t_{a_1}^\alpha\in \iM{N\bez U}=\gen{t_{a_2},t_{a_3}}. \]
In order to see that the above relation can not hold, consider $M$ as a subsurface of $\fal{M}=M_3$ as shown in Figure \ref{RN07}(ii). By Theorem \ref{tw:injectivity}, we can consider $\M{M}$ as a subgroup of $\M{\fal{M}}$ and now it is clear that the twist $t_c$ commutes with both $t_{a_2}$ and $t_{a_3}$, where $c$ is a circle indicated in Figure \ref{RN07}(ii). On the other, by Proposition \ref{prop:4.7}, $t_{a_1}^\alpha$ does not commute with $t_c$.
\end{ex}

\section{Virtually abelian geometric subgroups}

\begin{tw} \label{tw:VA}
Let $N$ be an essential subsurface of $M$. Then the geometric subgroup $\iM{N}$ is virtually
abelian if and only if $N$ is a disjoint union of Klein bottles with one hole, skirts, pantalons
and annuli.
\end{tw}
\begin{proof}
Suppose first that $N$ is a disjoint union of Klein bottles with one hole, skirts, pantalons and annuli.
Since the mapping class group of each of the listed surfaces is virtually abelian, $\iM{N}$ as a
homomorphic image of a virtually abelian group is virtually abelian (cf Corollary \ref{rem:prod} and Propositions \ref{prop:ex:virt}, \ref{prop:ex:semprod}).

Conversely, let $\iM{N}$ be virtually abelian and suppose that $N$ has a component, say $U$, which
is not a Klein bottle with a hole, skirt, pantalon nor annulus. Then there exists a generic two-sided
circle $a$ in $U$ which is not isotopic to a boundary component of $U$. By Proposition \ref{prop:gencirc},
there exists a generic two-sided circle $b$ in $U$ such that $I(a,b)>0$. On the other hand, since
$\iM{N}$ is virtually abelian, $i_*(t_a^\alpha)$ and $i_*(t_b^\beta)$ commute for some nonzero integers $\alpha$ and $\beta$. But $i_*(t_a^\alpha)$ and $i_*(t_b^\beta)$ are just $t_a^\alpha$ and $t_b^\beta$ treated as elements of $\M{M}$. Hence, by Proposition \ref{prop:4.7}, $I(a,b)=0$ which is a contradiction.
\end{proof}
\begin{wn}\label{wn:VirtAbel}
Let $N$ be a generic subsurface of $M\neq M_{-2}$ such that $\iM{N}$ is virtually abelian. Then there exist
two-sided generic circles $\lst{a}{n}$ in $N$ such that the subgroup
$\gen{t_{a_1},\ldots,t_{a_n}}$ is a finite index subgroup of $\iM{N}$ and is isomorphic to
$\zz^n$. Furthermore, the circles $\lst{a}{n}$ with the above properties are unique up to a
permutation or replacing one of the boundary components of an exterior
cylinder to $N$ with the second one.
\end{wn}
\begin{proof}
Define the circles $\lst{a}{n}$ in the following way:
\begin{itemize}
 \item take the meridian of each annulus of $N$,
 \item take the unique nonseparating two-sided circle on each component of $N$
     homeomorphic to a Klein bottle with one hole (cf Proposition \ref{prop:but:kl}),
 \item take all the boundary components of components of $N$ different from annuli
 \item for every cylinder $C$, exterior to $N$, remove from the above set of circles one of the boundary
 components of $C$, unless $M=M_1$.
\end{itemize}
By the previous theorem, $N$ is a disjoint union of Klein bottles with one hole, skirts,
pantalons and annuli, hence by Corollary \ref{rem:prod}, 
$\gen{t_{a_1},\ldots,t_{a_n}}$ is a finite index subgroup of $\iM{N}$ and by Proposition
\ref{prop:4.4}, \[\gen{t_{a_1},\ldots,t_{a_n}}\simeq \zz^n.\]

In order to show the uniqueness of circles $\lst{a}{n}$, observe that up to replacing one
of the boundary components of an exterior cylinder to $N$ with the second one, the set
$\{\lst{a}{n}\}$ contains all generic two-sided circles in $N$. Moreover, skipping any of the $a_i$'s leads to an infinite index subgroup of
$\gen{t_{a_1},\ldots,t_{a_n}}$.
\end{proof}


\section{Commensurability}

%

\begin{lem}\label{lem:tw:comm:1}
Let $N_0$ and $N_1$ be generic subsurfaces of $M$ which do not have isotopic components and do not have components homeomorphic to a Klein bottle with one hole. If the geometric subgroups $\iM{N_0}$ and $\iM{N_1}$ are commensurable then
\begin{enumerate}
 \item every boundary component of $N_i$ is isotopic in $M$ to a boundary component of $N_{1-i}$ for $i=0,1$;
 \item there exists a subsurface $S$ of $M$ such that $\kre{S\bez N_0}$ is isotopic to $N_1$ and for each boundary component $d$ of $S$ there exists a component of $N_0$ or $N_1$ which is an annulus with meridian isotopic to $d$;
 \item if $a$ is a generic two-sided circle in $N_i$ then $a$ is isotopic in $M$ to a circle disjoint from $N_{1-i}$ for $i=0,1$.
\end{enumerate}

\end{lem}

\begin{proof}
We first prove assertion (1). Using the symmetric role of $N_0$ and $N_1$, we can
concentrate on the case $i=0$.
Let $\lst{d}{n}$ denote the boundary components of $N_0$, and let
$\lst{a}{m}$ be the circles
which determine P-S decompositions of components of $N_0$ different from annuli. Define also $\lst{b}{k}$ to be
the union of circles defining P-S decompositions for components of $\kre{M\bez N_0}$ different from annuli, and the boundary
circles of $M$ not isotopic to any $d_i$ -- see Figure \ref{RN03}.
\begin{figure}[h]
\includegraphics{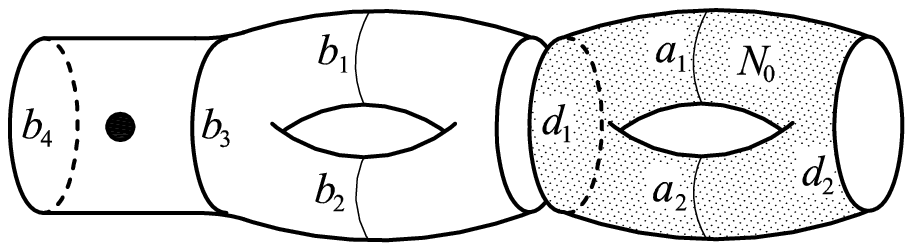}
\caption{Subsurface $N_0$ -- Lemma \ref{lem:tw:comm:1}.} \label{RN03}
\end{figure}
If $d$ is one of boundary components of $N_1$, then 
\begin{itemize}
 \item $I(d,a_i)=0$ for each $1\leq i \leq m$. In fact, by commensurability, $t_{a_i}^\alpha=i^{N_0}_*(t_{a_i}^\alpha)\in{i^{N_1}_*({\cal{M}}(N_1))}$ for some integer
$\alpha\neq 0$, and $t_{d}$ is central in ${i^{N_1}_*({\cal{M}}(N_1))}$. Hence by Proposition \ref{prop:4.7},
$I(d,a_i)=0$ (note that we use here the assumption that $N_1$, hence $d$, is generic).
 \item $I(d,b_j)=0$, $I(d,d_l)=0$ for each $1\leq j\leq k$, $1\leq l\leq n$.
 The reason is similar as before, $t_{d}^\alpha\in\iM{N_0} $ for some
 integer $\alpha\neq 0$ and each of $t_{b_j},t_{d_l}$ commutes with all of $\iM{N_0}$.
\end{itemize}
Therefore by Proposition \ref{prop_6_2}, we can assume that $d$ is
disjoint form each of $a_i,b_j,d_l$. Therefore $d$ lies entirely in one of the connected
components of $M\bez (\bigcup_{i=1}^{m} a_i\cup\bigcup_{j=1}^{k} b_j\cup\bigcup_{l=1}^{n} d_l)$.
Therefore $d$ is isotopic to one of $a_i,b_j,d_l$.

Now we are going show that $d$ can not be isotopic to $b_j$. Suppose on the contrary that
$d\simeq b_j$. By commensurability, $t_{d}^\alpha\in \iM{N_0}$ for some integer
$\alpha\neq 0$.
Hence by Lemma \ref{lem:comm:pow:tw}, $d\simeq b_j$ is
isotopic to a boundary component of $N_0$, which is a contradiction with the definition of $b_j$.

Next suppose that $d\simeq a_i$. By Proposition \ref{prop:gencirc}, there exists a generic two-sided circle $e$ in $N_0$, such that
$I(e,d)=I(e,a_i)>0$ (we use here the assumption that $N_0$ does not have a Klein bottle with one hole as a component). Now $t_e^{\alpha}\in{i_*({\cal{M}}(N_1))}$ for some integer $\alpha\neq
0$. Since $t_{d}$ is central in ${{\cal{M}}(N_1)}$, Proposition \ref{prop:4.7} implies that
$I(e,d)=0$, a contradiction.

Therefore, $d$ is isotopic to $d_l$ for some $l=1,\ldots,q$.
%

We now turn to the proof of (2). For $i\in \{0,1\}$ define \[{\cal{C}}_i=\{c_{i,1},\ldots,c_{i,n_i}\}\] to be the
smallest set of circles in $M$ with the following properties:
\begin{itemize}
\item for every component $U$ of $N_i$ which is not an annulus or which is an exterior cylinder to $N_{1-i}$, ${\cal{C}}_i$ contains boundary components of $U$;
\item for every component $U$ of $N_i$ which is an annulus and not an exterior cylinder to $N_{1-i}$, ${\cal{C}}_i$ contains this boundary component $d$ of $U$ which is closer to $N_{1-i}$, in the sense that after an isotopy of $U$ in $M$ which takes $d$ to a boundary component of $N_{1-i}$ (such an isotopy exists by statement (1)), the second boundary circle of $U$ is disjoint from $N_{1-i}$ (since $U$ is not an exterior cylinder of $N_{1-i}$, only one boundary component of $U$ can satisfy this condition).
\end{itemize}
By statement (1) and Proposition \ref{prop_6_10}, $n_0=n_1$ and there exists an isotopy of $M$ which takes ${\cal{C}}_0$ to
${\cal{C}}_1$. Therefore we can assume that ${\cal{C}}_0={\cal{C}}_1$.

 We claim that ${N_0}\cap {N_1}={\cal{C}}_0$. Suppose on the contrary that the interiors of $U_0$ and $U_1$ intersect, where $U_0$ and $U_1$ are some components of $N_0$ and $N_1$ respectively. If $U_1\not\podz U_0$ then some boundary component of $U_0$ intersects the interior of $U_1$, hence it is not in ${\cal{C}}_0$. Therefore $U_0$ is an annulus and we obtain a contradiction with the construction of ${\cal{C}}_0$. In fact, this contradicts the assumption that ${\cal{C}}_0$ contains the boundary component of $U_0$ which is ``closer'' to $N_1$.

Similarly we argue that $U_1\podz U_0$, hence $U_0=U_1$ which is a contradiction with the assumption that $N_0$ and $N_1$ do not have isotopic components.

Therefore ${N_0}\cap {N_1}={\cal{C}}_0$, and if we define $S=N_0\cup N_1$ then $N_1=\kre{S\bez N_0}$. Moreover, the boundary components of $S$ are exactly the boundary components of $N_0$ and $N_1$ which are not in ${\cal{C}}_0$, and by construction of ${\cal{C}}_0$, they are boundary components of annuli of $N_0$ and $N_1$.

Assertion (3) is an immediate consequence of (2).
\end{proof}

The following lemma can be thought as a partial converse to Proposition \ref{prop:ex:semprod} -- cf Example \ref{ex:semprod:opp}.

\begin{lem}\label{lem:comm:red:common}
Let $U_0,U_1,U$ be generic subsurfaces of $M$ such that $U_i\cap U=\emptyset$ and $N_i=U_i\cup U$ is a generic subsurface of $M$ for $i=0,1$ -- see Figure \ref{RN06}.
\begin{figure}[h]
\includegraphics{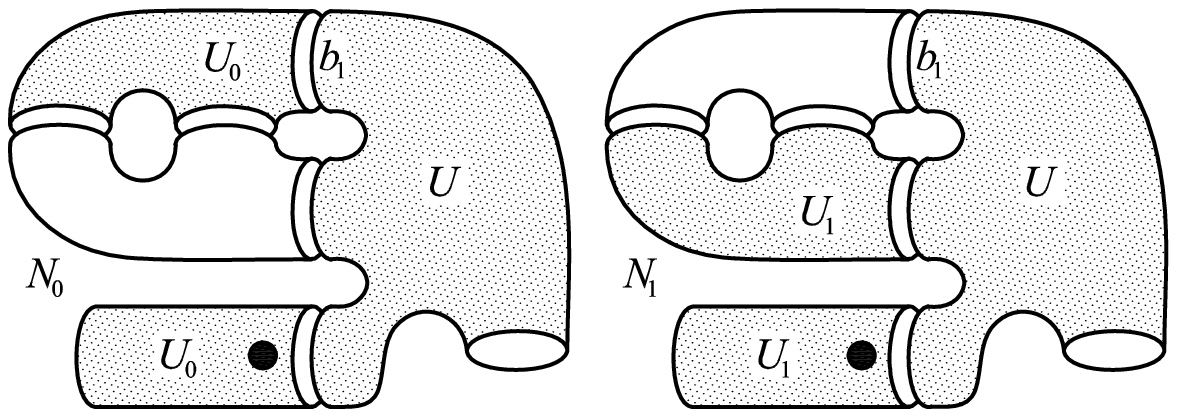}
\caption{Subsurfaces $N_0$ and $N_1$ -- Lemma \ref{lem:comm:red:common}.} \label{RN06}
\end{figure}
Let $\fal{U}_i$ be a surface obtained from $U_i$ by adding a regular neighbourhood of each boundary circle of $U$ which is isotopic to a boundary circle of $U_{1-i}$ but not isotopic to a boundary circle of $U_i$, for $i=0,1$.  Then $\fal{U}_0,\fal{U}_1$ are generic subsurfaces of $M$ and $\iM{N_0}$ is commensurable to $\iM{N_1}$ if and only if $\iM{\fal{U}_0}$ is commensurable to $\iM{\fal{U}_1}$.
\end{lem}
\begin{proof}
The ``if'' clause is an easy consequence of Corollary \ref{rem:prod} and Proposition \ref{prop:ex:semprod}, hence we concentrate on the ``only if'' part.

Assume that the groups $\iM{N_0}$ and $\iM{N_1}$ are commensurable.
Since $\iM{N_0}\cap\iM{N_1}$ has finite index in both $\iM{N_0}$ and $\iM{N_1}$, by Proposition \ref{prop:ex:cap},
\[\iM{N_0}\cap\iM{N_1}\cap \iM{\fal{U}_i}=\iM{N_{1-i}}\cap \iM{\fal{U}_i}\]
is of finite index in $\iM{\fal{U}_i}$ for $i=0,1$. Therefore it is enough to show that
\begin{equation}\label{eq:comm:lem}
\iM{N_{1-i}}\cap \iM{\fal{U}_i}=\iM{\fal{U}_{0}}\cap \iM{\fal{U}_1}\quad\text{for $i=0,1$}.
\end{equation}
Using the symmetric role of $N_0$ and $N_1$, we will restrict ourselves to the case $i=0$. Let $\lst{b}{k}$ be all the circles which are boundary circles of both $U$ and $U_0$, and which are not boundary circles of $U_1$ -- see Figure \ref{RN06}.

Applying Proposition \ref{prop:supp:disjoint2} with $N=N_{1}$, $U=U$, $K=\fal{U}_0$  and \[f\in \iM{N_{1}}\cap \iM{\fal{U}_0},\] we obtain
\[f=gt_{b_1}^{\beta_1}\cdots t_{b_k}^{\beta_k},\]
where $g\in \iM{U_{1}}$. Hence $f\in \iM{\fal{U}_1}$, which completes the proof of equality \eqref{eq:comm:lem}.
%
\end{proof}

\begin{tw} \label{tw:comm:va}
Let $N_0$ and $N_1$ be generic subsurfaces of $M$ such that no component of $N_0$ is isotopic to a
component of $N_1$. Then the geometric subgroups $\iM{N_0}$ and $\iM{N_1}$ are commensurable if and
only if they are virtually abelian with the same set of basic circles.
\end{tw}
\begin{proof}
The ``if'' clause follows from Corollary \ref{wn:VirtAbel}. Hence we concentrate on the ``only
if'' clause, that is assume that $\iM{N_0}$ and $\iM{N_1}$ are commensurable.

Observe that by Corollary \ref{wn:VirtAbel}, it is enough to show that both of these groups are
virtually abelian. In fact, if this is the case, then the set of basic circles for $\iM{N_0}$ is
also a set of basic circles for $\iM{N_1}$.

Our next claim is that it is enough to consider only the case when $N_0$ and $N_1$ do not have
connected components homeomorphic to a Klein bottle with one hole. In fact, suppose that $N_0$ has a
component $K$ homeomorphic to a Klein bottle with one hole. Let $\fal{N}_0$ be a surface obtained form
$N_0$ by removing $K$ and adding:
\begin{itemize}
 \item[-] regular neighbourhood of a nonseparating two-sided circle in $K$,
 \item[-] regular neighbourhood of a boundary component of $K$ if this boundary component is not
 isotopic to a boundary component of $N_0\bez K$.
\end{itemize}
Then $\fal{N}_0$ is a generic subsurface. Moreover, by Corollary \ref{rem:prod} and Propositions \ref{prop:but:kl} and \ref{prop:ex:semprod},  the geometric subgroups
$\iM{N_0}$ and $\iM{\fal{N}_0}$ are commensurable, hence it is enough to prove that both $\iM{\fal{N}_0}$ and
$\iM{N_1}$ are virtually abelian. Moreover, if $\fal{N}_0$ and $N_1$ have common (up to isotopy) connected
component $C$, then $C$ must be one of the annuli added to $N_0\bez K$. Hence by Lemma \ref{lem:comm:red:common}, it is enough to prove that $\iM{\fal{N}_0\bez C}$ and $\iM{N_1\bez C}$ are virtually
abelian. If $\fal{N}_0\bez C$ and $N_1\bez C$ still have a common connected component (i.e. the second of
the added annuli), we can remove it as before.

Therefore, by repeating the above procedure of removing Klein bottles with one hole, we can assume that neither
$N_0$ nor $N_1$ have components homeomorphic to a Klein bottle with one hole.


In order to finish the proof of the theorem, assume that some component $U$ of $N_0$
is not a pair of pants, skirt nor annulus. Then there exists a generic two-sided circle
$a$ in $U$ which is not isotopic to a boundary component. By statement (3) of Lemma \ref{lem:tw:comm:1}
and by Lemma \ref{lem:comm:pow:tw}, $a$ is isotopic to a
boundary circle of $N_1$, hence by statement (1) of Lemma \ref{lem:tw:comm:1}, $a$ is also a
boundary circle of $N_0$, which is a contradiction with the definition of $a$. Therefore, by Theorem
\ref{tw:VA}, $\iM{N_0}$ is virtually abelian. Clearly the same is true for
$\iM{N_1}$.
\end{proof}
The following example shows that the statement of Theorem \ref{tw:comm:va} would be significantly more complicated if we did not require $N_0$ and $N_1$ to be generic.
\begin{ex}\label{ex:notgeneric:subs}
Let $M$ be a torus with two punctures $z_0,z_1$ and let $N_i$ be a complement in $M$ of a small disk around $z_i$ for $i=0,1$. Since $\M{M}$ is generated by twists about the meridian and longitude of $M$, $\iM{N_0}=\iM{N_1}=\M{M}$ despite the fact that this group is not virtually abelian.
\end{ex}

\section{Commensurability - geometric interpretation}
Keeping in mind Theorem \ref{tw:comm:va} and Lemma \ref{lem:comm:red:common}, we obtain a general construction of generic subsurfaces $N_0$
and $N_1$ such that $\iM{N_0}$ and $\iM{N_1}$ are commensurable. In fact, start with two sets $N_0$ and
$N_1$ consisting of skirts, pantalons, annuli and Klein bottles with one hole. Then glue elements of $N_0$ to elements of $N_1$ along some of the boundary components, denote the obtained surface by $S$. In order to make the groups $\iM{N_0}$ and $\iM{N_1}$ commensurable we need to ensure that the boundary twists of $S$ and nonseparating two-sided circles in components of $S$ homeomorphic to Klein bottles with one hole are both in $\iM{N_0}$ and $\iM{N_1}$. The case of a generic nonseparating two-sided circle $a$ in a component $K$ of $N_i$ which is a Klein bottle with a hole, can be fixed by adding to $N_{1-i}$ either a regular neighbourhood of $a$ or a complement of a neighbourhood of $a$ in $K$, for $i=0,1$.
In order to fix the problem with boundary components of $S$ we can iterate the following technics:
\begin{itemize}
 \item we can add an arbitrary surface $U$, disjoint from $N_0\cup N_1$, to both $N_0$ and $N_1$; 
 \item suppose that $d$ is a boundary component of $S$ which is in ${N_i}$, then add to $N_{1-i}$ a regular neighbourhood of $d$, for $i=0,1$.
\end{itemize}
Finally, embed obtained surface $N_0\cup N_1$ in some surface $M$ in such a way that $N_0$ and $N_1$ are generic subsurfaces.

Our next goal is to prove that the described construction of $N_0$ and $N_1$ is as general as possible, i.e. that every pair of generic subsurfaces which lead to commensurable geometric subgroups can be constructed in that way. However, in order to simplify the formulation, we divide the statement into two steps (Theorems \ref{tw:GeomtrChar:pre} and \ref{tw:GeomtrChar} below). 
\begin{tw}\label{tw:GeomtrChar:pre}
Let $N_0$ and $N_1$ be generic subsurfaces of $M$ such that no component of $N_0$ is isotopic to a component of $N_1$. Assume also that no component of $N_0$ or $N_1$ is a Klein bottle with one hole. Then the geometric subgroups $\iM{N_0}$ and $\iM{N_1}$ are commensurable if and only if
\begin{enumerate}
\item each component of $N_0$ and $N_1$ is a skirt a pantalon or an annulus;
\item there exists a subsurface $S$ of $M$ such that  $\kre{S\bez N_0}$ is isotopic to $N_1$;
\item for each boundary component $d$ of $S$ there exists a component of $N_0$ or $N_1$ which is an annulus with meridian isotopic to $d$.
\end{enumerate}
\end{tw}
\begin{proof}
Conditions (1)--(3) clearly imply that $N_0$ and $N_1$ are virtually abelian with the same set of basic circles, which proves the ``if'' clause.

The ``only if'' clause is an immediate consequence of Theorem \ref{tw:comm:va} and assertion (2) of Lemma \ref{lem:tw:comm:1}.
\end{proof}
As an immediate corollary we obtain the following natural generalisation of Theorem 6.5 of \cite{RolPar}.
\begin{wn}\label{wn:comm:bez:annulus}
 Let $N_0$ and $N_1$ be generic subsurfaces of $M$ such that no component of $N_0$ is isotopic to a component of $N_1$.
 Assume also that neither $N_0$ nor $N_1$ has components homeomorphic to either a Klein bottle with one hole or an annulus.
 Then the geometric subgroups $\iM{N_0}$ and $\iM{N_1}$ are commensurable if and only if $M$ has no boundary, each component of $N_0$ and of $N_1$ is a pantaloon or a skirt and $\kre{M\bez N_1}$ is isotopic to $N_0$.\qed
\end{wn}
\begin{tw}\label{tw:GeomtrChar}
Let $N_0$ and $N_1$ be generic subsurfaces of $M$ such that the geometric subgroups $\iM{N_0}$ and $\iM{N_1}$ are commensurable. Let $\lst{U}{p}$ be the set of all common (up to isotopy) connected components of $N_0$ and $N_1$.
Let $K_{i,1},\ldots,K_{i,q_i}$ be the only connected components of $N_i'=N_i\bez\bigcup_{j=1}^{p}U_j$ which are homeomorphic to a Klein bottle with one hole, for $i=0,1$.
Then
there exists a component $L_{1-i,j}$ of $N_{1-i}'$ which is isotopic to either a regular neighbourhood of the nonseparating two-sided circle in $K_{i,j}$ or to the complement of this neighbourhood in $K_{i,j}$, for $j=1,\ldots q_{i}$ and $i=0,1$.

Moreover, for $i=0,1$, let $\fal{N}_i$ be a surface obtained from \[\widehat{N}_i=N_i'\bez\left( \bigcup_{j=1}^{q_i} K_{i,j} \cup \bigcup_{j=1}^{q_{1-i}} L_{i,j}\right)\]  as follows:
for each $U\in\{\lst{U}{p},K_{0,1},\dots,K_{0,q_0},K_{1,1},\dots,K_{1,q_1}\}$, and for each boundary component $d$ of $U$ which is isotopic in $M$ to a boundary component of $\widehat{N}_{1-i}$ and is not isotopic to a boundary component of $\widehat{N}_i$, add a regular neighbourhood of $d$ to $\widehat{N}_i$.
 Then $\fal{N}_0$ and $\fal{N}_1$ satisfy the assumptions of Theorems \ref{tw:comm:va} and \ref{tw:GeomtrChar:pre}. Moreover, $\iM{\fal{N}_0}$ is commensurable to $\iM{\fal{N}_1}$.
\end{tw}
\begin{proof}
By Lemma \ref{lem:comm:red:common}, up to adding some annuli, surfaces ${N_0'}$ and ${N_1'}$ satisfy the assumptions of Theorem \ref{tw:comm:va}. By Theorem \ref{tw:VA}, the components of ${N_0'}$ and ${N_1'}$ are Klein bottles with one hole, pantaloons, skirts or annuli. Moreover, by Corollary \ref{wn:VirtAbel}, the sets of basic circles coincide. Therefore, if $i=0,1$, $1\leq j\leq p$ and $a$ is the nonseparating two-sided circle in $K_{i,j}$, then $a$ is isotopic to a circle in some connected component $L_{1-i,j}$ of $N_{1-i}$. Therefore,  $L_{1-i,j}$ is isotopic to either a regular neighbourhood of $a$ or to its complement in $K_{i,j}$.

The fact that $\fal{N}_0$ and $\fal{N}_1$ do not have common connected components, and that they do not have Klein bottles with one hole as components is obvious. Moreover, by Lemma \ref{lem:comm:red:common}, $\iM{\fal{N}_0}$ is commensurable to $\iM{\fal{N}_1}$.
\end{proof}
\begin{wn}
Let $N_0$ and $N_1$ be generic subsurfaces of $M$ which do not have skirts, pantalons nor annuli as connected components. Then the geometric subgroups $\iM{N_0}$ and $\iM{N_1}$ are commensurable if and only if $N_0$ is isotopic to $N_1$.\qed
\end{wn}
%


\section{Commensurator of a geometric subgroup}

Our next goal is to describe the commensurator ${\rm{Comm}}(\iM{N})$ of the geometric subgroup corresponding to an arbitrary generic subsurface $N$. The first guess could be that this commensurator should be close to the \emph{stabiliser of $N$}
\[{\rm Stab}(N)=\{f\in \M{M}\st \text{$f(N)$ is isotopic to $N$}\} \]
 In fact it is known that this is the case for $M$ orientable, and $N$ injective and connected -- cf Theorem 7.1 of \cite{RolPar}. However, as is shown in the following example, the index \[ [{\rm{Comm}}(\iM{N}):{\rm{Stab}}(N)]\] can be arbitrary large.
\begin{ex}\label{ex:index:comm:large}
For $n\geq 3$, let $M$ be a sphere with $3n$ punctures, embedded in $\rr^3$ in a rotationally symmetric manner indicated in Figure~\ref{RN04} (this figure shows the case $n=3$, in order to imagine the general case just think of $n$ ``branches'' instead of $3$).
\begin{figure}[h]
\includegraphics{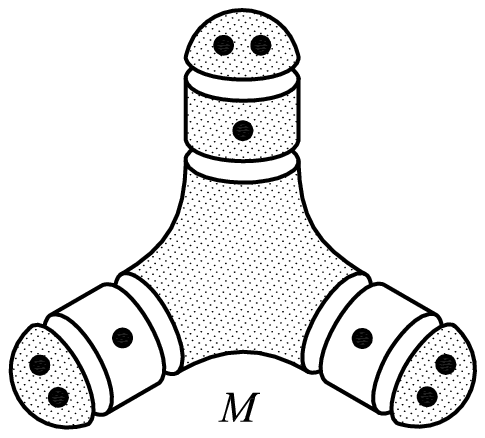}
\caption{Sphere with $3n$ punctures -- Example \ref{ex:index:comm:large}.} \label{RN04}
\end{figure}
The same figure indicates a generic subsurface $N$ (the shaded region), it consists of $n$ ``outer'' doubly-punctured disks, ``the core''  of $M$ (which is a sphere with $n$ holes) and one punctured annulus. Now it is clear that the natural rotations of $M$ provide $n$ elements of ${\rm{Comm}}(\iM{N})$ (even of the normaliser of $\iM{N}$) that represent different cosets of ${\rm{Stab}}(N)$. Hence
\[[{\rm{Comm}}(\iM{N}):{\rm{Stab}}(N)]\geq n\]
\end{ex}
Despite the above example, we can still provide a descent characterisation of ${\rm{Comm}}(\iM{N})$ but we need to ``redefine'' the stabiliser of $N$.
\begin{df}
Let $\lst{U}{p}$ be all connected components of a generic subsurface $N$ which are not annuli, skirts, pantalons nor Klein bottles with one hole. Denote also by ${\cal{C}}$ the set of isotopy classes of the union of basic circles for the subsurface $N\bez\bigcup_{i=1}^p U_i$ and the boundary circles of $\bigcup_{i=1}^p U_i$. Define ${\rm{Stab}}^*(N)$ to be the subgroup of $\M{M}$ consisting of these classes of diffeomorphisms $\map{f}{M}{M}$ for which $f\left(\bigcup_{i=1}^p U_i\right)$ is isotopic to $\bigcup_{i=1}^p U_i$, and $f({\cal{C}})$ is isotopic to ${\cal{C}}$.
\end{df}
\begin{tw} \label{tw:decr:comm1}
Let $N$ be a generic subsurface of $M$. Then
\[{\rm{Comm}}(\iM{N})={\rm{Stab}}^*(N). \]
\end{tw}
\begin{proof}
The inclusion ${\rm{Stab}}^*(N)\podz{\rm{Comm}}(\iM{N})$ is obvious, hence we concentrate on the second one. Let $f\in {\rm{Comm}}({\iM{N}})$, that is $\iM{N}$ and $\iM{f(N)}$ are commensurable. By Theorem \ref{tw:GeomtrChar}, the subsurfaces $\fal{N}$ and $\fal{f(N)}$ (constructed as in the statement of that theorem) have commensurable geometric subgroups, no common components and no components homeomorphic to a Klein bottle with one hole. Hence by Theorems \ref{tw:comm:va} and \ref{tw:VA}, each component of $\fal{N}$ and of $\fal{f(N)}$ is an annulus, a skirt or a pantaloon. Moreover, the sets of basic circles for these surfaces must coincide which easily leads to the conclusion that $f\in {\rm{Stab}}^*(N)$.
\end{proof}

\begin{tw}\label{tw:decr:comm2}
Let $N$ be an injective subsurface of $M$ such that no component of $N$ is a Klein bottle with
one hole or an annulus.
 Then
\begin{enumerate}
\item ${\rm Comm}(\iM{N})={\rm Stab}(N)\rtimes \zz_2$ if $M$ is closed, $\iM{N}$ is virtually abelian and there exists a diffeomorphism \[\map{\sig}{N}{\kre{M\bez N}}\] such that $\sig\in \M{M}$; 
\item ${\rm Comm}(\iM{N})={\rm Stab}(N)$ otherwise.
\end{enumerate}
\end{tw}
\begin{proof}
The inclusion ${\rm Stab}(N)\podz{\rm Comm}(\iM{N})$ is obvious, hence it is enough to prove that
if there exists an element \[\sig \in {\rm Comm}(\iM{N})\bez {\rm Stab}(N)\] then ${\rm
Comm}(\iM{N})={\rm Stab}(N)\rtimes \zz_2$ and $N$ is as described in (1).

Let ${\cal U}=\{\lst{U}{p}\}$ be the set of all components of $N$ such that $\sig \in {\rm Stab}
(U_i)$ for $i=1,\ldots, p$. By Lemma \ref{lem:comm:red:common}, $N'=N\bez {\cal U}$ and $\sig(N')=\sig(N)\bez
{\cal U}$ have commensurable geometric subgroups (we use here injectivity of $N$).
Hence $N'$ and $\sig(N')$ satisfy the assumptions of Corollary \ref{wn:comm:bez:annulus},  which
implies that $M$ is closed, $\sig(N')$ is isotopic to $\kre{M\bez N'}$ and $N'$ is virtually abelian.
Therefore, ${\cal U}=\emptyset$, $N=N'$ and we have an exact sequence
\[\begin{CD}1@>>>{\rm Stab}(N)@>>> {\rm Comm}(\iM{N})@>\pi>> \zz_2@>>>1\end{CD}\]
where $\zz_2=\{1,-1 \}$ and $\pi(h)=-1$ iff $h(N)$ is isotopic to $\kre{M\bez N}$ for $h\in {\rm Comm}(\iM{N})$.

It remains to show that the above sequence splits. In order prove this, embed $N$ in
\[\rr^4=\{(x_1,x_2,x_3,x_4)\st x_1,x_2,x_3,x_4\in\rr\}\]
in such a way that:
 \begin{itemize}
 \item the interior of $N$ is contained in the set $x_4< 0$,
 \item the boundary of $N$ is contained in the plane $x_3=x_4=0$,
 \item each boundary component of $N$ is a metric circle with the center on the $x_1$ axis.
 \end{itemize}
Now if $\map{\sig}{\rr^4}{\rr^4}$ is the half turn about the $x_1$-axis, that is
\[\sig(x_1,x_2,x_3,x_4)=(x_1,-x_2,-x_3,-x_4), \]
 then $N\cup_{\sig} \sig(N)$ is a
model for $M$ in which $-1\mapsto \sig$ provides a section
\[\map{s}{\zz_2}{{\rm Comm}(\iM{N})}\]
of $\pi$ defined above.
\end{proof}
\begin{uw}
One can easily see that the semi-direct product in part (1) of the above theorem is a direct
product if and only if each component of $N$ is a pantalon of type II or III, or else a skirt of
type II.
\end{uw} 
\section{Irreducible representations}\label{sec:rep}
As we indicated in the introduction, the mapping class group $\M{M}$ acts on the complex of curves $C(M)$ and this action produces a very interesting and important family of subgroups of $\M{M}$, namely the family of stabilisers of simplexes in $C(S)$. Since such a stabiliser ${\rm Stab}(\sigma)$ contains a geometric subgroup $\iM{M_\sigma}$ as a subgroup of finite index (cf Example \ref{ex:complex1}), we can apply the results from previous sections to the study of these stabilisers.
In particular we will show that under some natural assumptions,
\begin{equation}
{\rm{Stab}}(\sig_0)\text{ is commensurable to }{\rm{Stab}}(\sig_1) \Longleftrightarrow \sig_0=\sig_1.\label{eq:com:Stab}
\end{equation}
Then using the results from \cite{Burger}, we will draw some interesting conclusions concerning irreducible unitary representations of $\M{M}$.
In the orientable case, equivalence \eqref{eq:com:Stab} was proved in \cite{ParRep} by means of the ``large action'' of the mapping class group on $C(S)$ (in the sense of \cite{Burger}). However, in the nonorientable case there are new phenomena we need to deal with.

The first one is that keeping in mind Example \ref{ex:notgeneric:subs}, we do not want to deal with subsurfaces which are not generic, hence we will restrict ourselves to the complex $C_0(M)$ of two-sided circles in $M$ (cf Example \ref{ex:complex1}). This is a full subcomplex of $C(M)$ in the sense that if the set of vertices in $C_0(M)$ is a simplex in $C(M)$ then it is a simplex in $C_0(M)$.

The second phenomenon is illustrated by the following example.
\begin{ex}\label{ex:simpl:comm:stab}
Let $L$ be a surface with two boundary components $c_1$ and $c_2$ which is not an annulus. Construct a surface $M$ by gluing the boundary circle of a Klein bottle $K_i$ with one hole to $c_i$, for $i=1,2$, and let $a_1,a_2$ be the nonseparating two-sided circles in $K_1$ and $K_2$ respectively -- see Figure \ref{RN05}.
\begin{figure}[h]
\includegraphics{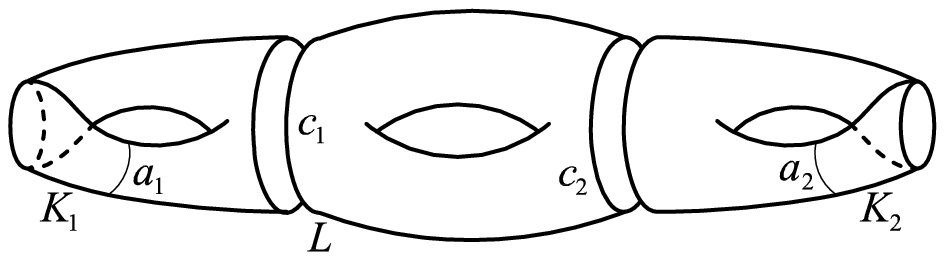}
\caption{Simplexes of $C_0(M)$ with commensurable stabilisers -- Example \ref{ex:simpl:comm:stab}.} \label{RN05}
\end{figure}
If we define the simplexes $\sig_0=\{c_1,c_2,a_1\}$ and $\sig_1=\{c_1,c_2,a_1,a_2\}$ then by assertion (1) of Proposition \ref{prop:but:kl}, ${\rm Stab}(\sig_0)\podz{\rm Stab}(\sig_1)$. Moreover, ${\rm Stab}(\sig_0)$ is the kernel of a homomorphism \[\map{\Phi}{{\rm Stab}(\sig_1)}{\zz_2}\] defined by $\Phi(f)=-1$ iff $f$ interchanges $c_1$ and $c_2$ for $f\in {\rm Stab}(\sig_1)$. Therefore $[{\rm Stab}(\sig_1):{\rm Stab}(\sig_0)]=2$, hence these stabilisers are commensurable despite the fact that $\sig_0\neq\sig_1$.
\end{ex}
The above example motivates the following definition.
\begin{df}
A simplex $\sig=\{a_0,\ldots,a_n\}$ of the complex $C_0(M)$ is \emph{reduced} if for each $i=0,\ldots, n$, no component of $\kre{M\bez a_i}$ is a Klein bottle with one hole. The \emph{reduced complex of curves}, denoted by $C^{red}_0(M)$ is a subcomplex of $C_0(M)$ consisting of reduced simplexes.
\end{df}
As is shown in the next proposition, the complex $C_0^{red}(M)$ is a very natural subcomplex of $C_0(M)$.
\begin{prop}
Suppose that $M\neq M_{-4}$. Then the realisation of the reduced complex of curves $|C_0^{red}(M)|$ is a strong deformation retract of the realisation $|C_0(M)|$ of the complex of two-sided curves.
\end{prop}
\begin{proof}
For any vertex $v$ of $C_0(M)$ define
\[\Phi(v)=\begin{cases}
v&\text{if no component of $\kre{M\bez v}$ is a Klein bottle with one hole;}\\
c&\text{if one of the components of $\kre{M\bez v}$ is a Klein bottle $K$}\\
&\text{with one hole, and $c$ is the unique}\\
&\text{nonseparating, two-sided circle in $K$.}
\end{cases} \]
Note that by our assumption $M\neq M_{-4}$, there is no ambiguity in the above definition. We claim that $\Phi$ can be extended to a simplicial map $\map{\Phi}{C_0(M)}{C_0^{red}(M)}$. In order to show this, we need to check that if $\sig=\{c_0,\dots,c_n\}$ is a simplex in $C_0(M)$ then $\{\Phi(c_0),\dots,\Phi(c_n)\}$ is a simplex in $C_0^{red}(M)$ (possibly of smaller dimension). Hence it is enough to show that $I(\Phi(c_i),\Phi(c_j))=0$ for $i\neq j$. This is trivial for circles which do not cut off a Klein bottle with one hole, so we can assume that one component $K$ of $\kre{M\bez c_i}$ is a Klein bottle with one hole. Then $c_j$ is either the unique nonseparating two-sided circle in $K$ or $c_j$ is a two-sided circle in $M\bez K$. In the first case $\Phi(c_i)=\Phi(c_j)$ and in the second case, using Lemma \ref{lem:zaw:but:kle}, it is straightforward to check that $\Phi(c_j)$ must be contained in $M\bez K$, hence $\Phi(c_j)$ is disjoint from $\Phi(c_i)$.

It remains to show that $\Phi$ is a strong deformation retraction. By the basic properties of weak topology (cf Chapter 3 of \cite{Spanier_book}), it is enough to consider the restriction of $\Phi$ to an arbitrary simplex $\sig\in C_0(M)$. By Proposition \ref{prop:but:kl} and Lemma \ref{lem:zaw:but:kle}, $\sig'=\sig\cup \Phi(\sig)$ is also a simplex in $C_0(M)$ and $\Phi$ restricted to $\sig'$ is just a simplicial projection onto the face of $\sig'$ spanned by vertices that are in $C_0^{red}$.
\end{proof}

\begin{prop} \label{prop:comm:symplexes:main}
Let $\sig_0,\sig_1\in { C_0^{red}}(M)$ be reduced simplexes such that ${\rm Stab}(\sig_0)$ and ${\rm Stab}(\sig_1)$ are commensurable. Then
$\sig_0=\sig_1$.
\end{prop}
\begin{proof}
It is an easy observation that for any $\sig\in C_0(M)$, ${\iM{M_\sig}}$ is a finite index subgroup of ${\rm Stab}(\sig)$, where
$M_\sig$ is defined as in Example \ref{ex:complex1}. Hence the geometric subgroups $\iM{M_{\sig_0}}$ and $\iM{M_{\sig_1}}$ are commensurable and neither $M_{\sig_0}$ nor $M_{\sig_1}$ contains components homeomorphic to an annulus or a Klein bottle with one hole. If $\lst{U}{p}$ are all (up to isotopy) common connected components of $M_{\sig_0}$ and $M_{\sig_1}$, then by Theorem \ref{tw:GeomtrChar}, ${M}_{\sig_0}'=M_{\sig_0}\bez\bigcup_{i=1}^p U_i$ and ${M}_{\sig_1}'=M_{\sig_1}\bez\bigcup_{i=1}^p U_i$ satisfy the assumptions of Theorem \ref{tw:comm:va}. Therefore $\iM{{M}_{\sig_0}'}$ and $\iM{{M}_{\sig_1}'}$ are virtually abelian with the same set of basic circles, which easily leads to the conclusion that $\sig_0=\sig_1$.
\end{proof}

As we have already indicated, the above proposition has very interesting consequences in terms of unitary irreducible representations of mapping class groups. In order to state the result, we need to recall the appropriate terminology and notation -- see \cite{Burger}.

For a countable discrete group $G$,
let $\widehat{G}$ be the \emph{unitary dual} of $G$, that is the set of equivalence classes of irreducible unitary representations of $G$. By $\widehat{G}^{fd}$ we denote the subspace of $\widehat{G}$ of equivalence classes of finite dimensional representations. Among the basic properties of irreducible representations is the following proposition (for a more complete statement and references see \cite{Burger}).
\begin{prop}[Mackey \cite{Mackey}]\label{tw:Mackey}
Let $\tau_i\in\widehat{H_i}^{fd}$ be a finite dimensional irreducible unitary representation of a subgroup $H_i$ of a countable and discrete group $G$ for $i=0,1$. Assume also that ${\rm Comm}(H_i)=H_i$ for $i=0,1$.  Then
\begin{enumerate}
 \item the induced representation $Ind^G_{H_i}(\tau_i)$ is irreducible for $i=0,1$, hence we have a well defined injective map
 \[\map{{\rm Ind}_{H_i}^G}{\widehat{H_i}^{fd}}{\widehat{G}};\]
 \item if in addition $H_0$ and $H_1$ are not conjugate in $G$, then $Ind^G_{H_0}(\tau_0)$ and $Ind^G_{H_1}(\tau_1)$ are not equivalent, hence we have a well defined injective map
     \[{\widehat{H_0}^{fd}}\sqcup{\widehat{H_1}^{fd}}\hookrightarrow{\widehat{G}}.\]
\end{enumerate}\qed
\end{prop}
A careful reader may noticed that in the original \cite{Burger} statement of (2) there is an assumption that $H_0$ and $H_1$ are not quasi-conjugate. However it is not hard to check that under our assumption ${\rm Comm}(H_i)=H_i$, $H_0$ and $H_1$ are quasi-conjugate if and only if they are conjugate.

As is shown in \cite{Burger}, there is an efficient way of constructing large families of subgroups of $G$ satisfying the assumptions of the above proposition by means of so called \emph{N.C.S. actions} of $G$.
\begin{df}[Burger, de la Harpe \cite{Burger}]
Let $G$ be a countable discrete group acting on a space $X$. The action of $G$ on $X$ is an \emph{action with noncommensurable stabilisers} (\emph{N.C.S. action} in short) if different points of $X$ have noncommensurable stabilisers.
\end{df}
\begin{prop}[Burger, de la Harpe \cite{Burger}]
Let $G\times X\to X$ be a N.C.S. action of a countable discrete group $G$. Then unitary induction provides a well defined injective map
\[\coprod_{x\in X/G} \widehat{{\rm Stab}(x)^{fd}}\hookrightarrow \widehat{G}\]
where $X/G$ is the orbit space.\qed
\end{prop}
Combining the above proposition with Proposition \ref{prop:comm:symplexes:main} leads to the following corollary.
\begin{wn}\label{cor:induc:repr}
The action of the mapping class group $\M{M}$ on the reduced complex of curves ${\cal C}^{red}_0(M)$ has
noncommensurable stabilisers in the sense of \cite{Burger}. Therefore unitary induction provides a
well defined injective map
\[\coprod_{\sig\in {\cal C}^{red}_0(M)/\M{M}} \widehat{{\rm Stab}(\sig)^{fd}}\hookrightarrow \widehat{\M{M}}\]\qed
\end{wn}

As we indicated before, the above corollary was proved in \cite{ParRep} for an orientable $M$.

\bibliographystyle{abbrv}
\bibliography{mybib}

\begin{thebibliography}{10}

\bibitem{Bir-Punct}
J.~S. Birman.
\newblock Mapping class groups and their relationship to braid group.
\newblock {\em Comm. Pure Appl. Math.}, 22:213--238, 1969.

\bibitem{Burger}
M.~Burger and P.~de~la {H}arpe.
\newblock Constructing irreducible representations of discrete groups.
\newblock {\em Proc. Indian Acad. Sci.}, 107:223--235, 1997.

\bibitem{Chil}
D.~R.~J. Chillingworth.
\newblock A finite set of generators for the homeotopy group of a
  non--orientable surface.
\newblock {\em Math. Proc. Cambridge. Philos. Soc.}, 65:409--430, 1969.

\bibitem{Epstein}
D.~B.~A. {E}pstein.
\newblock Curves on 2--manifolds and isotopies.
\newblock {\em Acta Math.}, 115:83--107, 1966.

\bibitem{Kork-non1}
M.~{K}orkmaz.
\newblock First homology group of mapping class groups of nonorientable
  surfaces.
\newblock {\em Math. Proc. Cambridge. Philos. Soc.}, 123(3):487--499, 1998.

\bibitem{Kork-non}
M.~{K}orkmaz.
\newblock Mapping class groups of nonorientable surfaces.
\newblock {\em Geom. Dedicata}, 89:109--133, 2002.

\bibitem{ParLab}
C.~{L}abru{\`{e}}re and {L}uis {P}aris.
\newblock Presentations for the puctured mapping class groups in terms of
  {A}rtin groups.
\newblock {\em Algebr. Geom. Topol.}, 1:73--114, 2001.

\bibitem{Lick1}
W.~B.~R. Lickorish.
\newblock A representation of orientable combinatorial 3-manifolds.
\newblock {\em Ann. of Math.}, 76:531--540, 1962.

\bibitem{Lick3}
W.~B.~R. Lickorish.
\newblock Homeomorphisms of non--orientable two--manifolds.
\newblock {\em Math. Proc. Cambridge. Philos. Soc.}, 59:307--317, 1963.

\bibitem{Lick2}
W.~B.~R. Lickorish.
\newblock A finite set of generators for the homeotopy group of a 2-manifold.
\newblock {\em Math. Proc. Cambridge. Philos. Soc.}, 60:769--778, 1964.

\bibitem{Mackey}
G.~W. Mackey.
\newblock {\em The theory of unitary group representations}.
\newblock Chicago Lectures in Mathematics. The University of Chicago Press,
  1976.

\bibitem{ParRep}
L.~{P}aris.
\newblock Actions and irreducible representations of the mapping class group.
\newblock {\em Math. Ann.}, 322:301--315, 2002.

\bibitem{RolPar}
L.~{P}aris and D.~{R}olfsen.
\newblock Geometric subgroups of mapping class groups.
\newblock {\em J. Reine Angew. Math.}, 521:47--83, 2000.

\bibitem{Spanier_book}
E.~H. Spanier.
\newblock {\em Algebraic topology}.
\newblock McGraw-Hill Book Co., 1966.

\bibitem{Stukow_SurBg}
M.~Stukow.
\newblock Mapping class groups of nonorientable surfaces with boundary.
\newblock Preprint.

\bibitem{Stukow_twist}
M.~Stukow.
\newblock Dehn twists on nonorientable surfaces.
\newblock {\em Fund. Math.}, 189:117--147, 2006.

\end{thebibliography}

\end{document}